\documentclass[12pt]{amsart}
 
\usepackage{amsmath,amsthm,amsfonts,amssymb}

\newcommand {\Z} {{\mathbb  Z}}
\newcommand {\R} {{\mathbb  R}}
\newcommand {\N} {{\mathbb  N}}

\newcommand {\ep} {\epsilon}
\newcommand {\tn} {\tfrac{2^{n+2}}{3}\pi}
\newcommand {\two} {\tfrac{\epsilon}{2^{n+2}}}
\newcommand {\thr} {\tfrac{\epsilon}{2^{n+3}}}

\newtheorem{theorem}{Theorem}

\newtheorem{proposition}{Proposition}

\allowdisplaybreaks[4]

\begin{document}
\title[Dimension Functions of Band-limited Wavelets]
{Estimation of Dimension Functions of Band-limited Wavelets}
\author{Biswaranjan Behera}
\address{Statistics and Mathematics Unit,
Indian Statistical Institute,
203, B. T. Road,
Calcutta-700108, India}
\email{br\_behera@yahoo.com}
\date{December 4, 2001\\ \hspace*{3.3mm} 2000 {\it Mathematics Subject Classification}. 42C40.}
%\subjclass[2000]{42C40}
\keywords{wavelet; band-limited; dimension function}

\begin{abstract}
The dimension function $D_\psi$ of a band-limited wavelet $\psi$ is bounded by $n$ if $\hat\psi$ is supported in 
$[-\tn,\tn]$. For each $n\in\N$ and for each $\ep$, $0<\ep<\delta=\delta(n)$, we construct a
wavelet $\psi$ with supp~$\hat\psi\subseteq$
$[-\tn,\tn+\ep]$ such that $D_\psi>n$ on a set of positive
measure, which proves that $[-\tn,\tn]$ is the largest symmetric
interval for estimating the dimension function by $n$. This construction also provides a family of (uncountably many) wavelet sets each consisting of infinite number of intervals.
\end{abstract}

\maketitle

\section{Introduction}
A wavelet is a function $\psi\in L^2(\R)$ such that the system $\{\psi_{j,k}=2^{j/2}\psi(2^j\cdot-k):j,k\in\Z\}$ forms an orthonormal basis for $L^2(\R)$. Given a wavelet $\psi$ of $L^2(\R)$, there is an associated function $D_\psi$, called the dimension function of $\psi$,
defined by
\begin{equation}\label{dpsi}
  D_\psi(\xi)=\sum\limits_{j\geq 1}\sum\limits_{k\in\Z}
  |\hat\psi(2^j(\xi+2k\pi))|^2.
\end{equation}

A simple periodization argument shows that  
$\int_0^{2\pi} D_\psi(\xi) d\xi= 2\pi\|\psi\|_2^2$, 
if $\psi\in L^2(\R)$. So the function $D_\psi$ is well defined and is finite a.e. Observe that $D_\psi$ is $2\pi$-periodic. 
P.\ G.\ Lemari\'e~\cite{Lem1, Lem2} used this function to show that certain wavelets are associated with a multiresolution analysis (MRA) of $L^2(\R)$. P.\ Auscher~\cite{Aus} proved that if $\psi$ is a wavelet, then the function $D_\psi$ is the dimension of certain closed subspaces of the sequence space $l^2(\Z)$ (hence the name dimension function, a term coined by Guido Weiss). This result in particular proves that $D_\psi$ is integer valued a.e. G.\ Gripenberg~\cite{Gri} 
and X.\ Wang~\cite{Wan}, independently, characterized all wavelets of $L^2(\R)$ associated with an MRA.
This well known characterization states that a wavelet $\psi$ of
$L^2(\R)$ is associated with an MRA if and only if $D_\psi=1$ a.e. The
article~\cite{BRS} contains a characterization of all dimension
functions.

A function is said to be band-limited if its Fourier transform is
compactly supported. It is easy to see that the dimension function
of a band-limited wavelet is bounded.

\begin{proposition}
\label{prop:one}
Let $n\in \N$. If $\psi$ is a wavelet such that supp~$\hat\psi\subset[-2n\pi,2n\pi]$, then $D_\psi\leq n$ a.e.
\end{proposition}

The above proposition is not optimal. For example, $D_\psi$ is still bounded by 1 for wavelets $\psi$ such that supp~$\hat\psi\subseteq [-\tfrac{8}{3}\pi, \tfrac{8}{3}\pi]$, which is proved in~\cite{HW} (see section~3.4). The authors of~\cite{BGRW} constructed an example of a wavelet $\psi$ with supp~$\hat\psi\subseteq$
$[-\tfrac{8}{3}\pi,\tfrac{8}{3}\pi+\epsilon],
0<\epsilon<\tfrac{2}{3}\pi$, such that $D_\psi\geq 2$ a.e.\ on a set
of positive measure, which shows that $[-\tfrac{8}{3}\pi,
\tfrac{8}{3}\pi]$ is the largest symmetric interval for estimating
the dimension function by 1. A natural question to ask is whether
there are optimal symmetric intervals to estimate the dimension
function by $n,\ n\geq 2$. The following theorem sheds light
to the above question.

\begin{theorem}\label{thm:opt}
Let $n\in \N$. If $\psi$ is a wavelet such that $\hat\psi$ is supported in $[-\tn,\tn]$, then $D_\psi\leq n$ a.e.
\end{theorem}

This result was also proved by Z.\ Rzeszotnik and D.\ Speegle in an unpublished article.
%(On the possible behaviour of wavelets in the frequency domain). 
They also proved that for every positive integer $n$ and every $\ep$, $0<\ep<\delta(n)$, there exists an MSF wavelet $\psi$ such that 
supp~$\hat\psi\subset[-\tn,\tn+\ep]$ and
$\|D_\psi\|_\infty>n$. This shows that $[-\tn,\tn]$ is the
optimal symmetric interval for estimating the dimension function
by $n$. We thank Professor Guido Weiss for kindly providing the above information to us. The purpose of this article is to construct such a wavelet explicitly. We shall prove the following theorem in a constructive manner.

\begin{theorem}
\label{thm:main}
For each $n\in\N$ and $0<\ep<\delta=\delta(n)$, there exists a wavelet $\psi$ such that supp~$\hat\psi\subseteq[-\tn,\tn+\ep]$ and $\|D_\psi\|_\infty>n$.
\end{theorem}

A wavelet $\psi$ of $L^2(\R)$ is said to be a {\it minimally
supported frequency} (MSF) wavelet if $|\hat\psi|$ is the characteristic function of some measurable subset $K$ of $\R$. The associated set $K$ is
called a {\it wavelet set}. A simple characterization of such sets is the following (see~\cite{HW} for a proof):

\begin{quote}
{\it A set $K\subset\R$ is a wavelet set if and only if both the
collections $\{K+2k\pi:k\in\Z\}$ and $\{2^jK:j\in\Z\}$ are
partitions of $\R$}.
\end{quote}

It is not always easy to construct wavelet sets satisfying desired properties. The concepts of translation and dilation equivalence of subsets of $\R$ are useful for this purpose.
% to construct wavelet sets. 
A set $A\subset\R$ is said to be $2\pi$-{\it translation equivalent} to a set $B\subset\R$ if there exists a partition $\{A_n:n\in\Z\}$ of $A$ such that $\{B_n\equiv A_n+2n\pi:n\in\Z\}$ is a partition of $B$. Similarly, $A$ is said to be $2$-{\it dilation equivalent} to $B$ if there exists another partition $\{A_n':n\in\Z\}$ of $A$ such that 
$\{B_n'\equiv 2^n A_n':n\in\Z\}$ is a partition of $B$.

In view of the characterization of wavelet sets stated above, it is now clear that a subset $K$ of $\R$ is a wavelet set if and only if $K$ is $2\pi$-translation equivalent to some interval of
length $2\pi$, $K\cap(0,\infty)$ is 2-dilation equivalent to
$[a,2a]$ for some $a>0$, and $K\cap(-\infty,0)$ is 2-dilation
equivalent to $[-2b,-b]$ for some $b>0$.
%%==============================
\section{Proofs of the theorems}
%\noindent{\it Proof of Proposition {\rm\ref{prop:one}}.} 
\proof[Proof of Proposition {\rm\ref{prop:one}}] 
Let $F(\xi)=\sum_{j\geq 1}|\hat\psi(2^j\xi)|^2$.
The condition on the support of $\hat\psi$ implies that 
 supp~$F\subset [-n\pi,n\pi]$. Since the equality,
\[
\sum_{j\in\Z}|\hat\psi (2^j\xi)|^2 = 1 \quad {\rm for~a.e.}\ \xi\in\R,
\]
is satisfied by every wavelet $\psi$, we have $F\leq 1$.  Therefore, we get $F\leq\chi_{[-n\pi ,n\pi]}$. This implies that
\[
D_{\psi}(\xi) = \sum_{k\in\Z}F(\xi + 2k\pi)
\leq\sum_{k\in\Z}\chi_{[-n\pi ,n\pi]}(\xi + 2k\pi) = n,
\]
 which proves the proposition.
\qed

\proof[Proof of Theorem {\rm\ref{thm:opt}}] 
Since the function $D_\psi$ is $2\pi$-periodic, it is enough to 
prove that if $\psi$ satisfies the hypothesis, then
$D_\psi(\xi)\leq n$ for $\xi\in[-\pi,\pi]$. For
 $\xi\in [-\pi,\pi]$, we have
 $(2k-1)\pi\leq \xi+2k\pi \leq (2k+1)\pi$ for all $k\in\Z$.

\smallskip 

(i) {$j=n$. If $k\geq 2$, then  $2^j(\xi + 2k\pi) \geq 2^j(2k-1)\pi=  2^n(2k-1)\pi\geq 3\cdot 2^n\pi \geq\tfrac{2^{n+2}}{3}\pi$. Similarly, if $k\leq -2$, then  $2^n(\xi +2k\pi)\leq -\tfrac{2^{n+2}}{3}\pi$. Hence, for $j=n,$ the only non-zero terms contributing to $D_{\psi}$ are for $k=-1,0,1$.

\smallskip

 (ii) $j\geq n+1$. If $k\geq 1$, then
$2^j(\xi + 2k\pi)\geq 2^j(2k-1)\pi\geq 2^{n+1}(2k-1)\pi\geq
2^{n+1}\pi \geq \tfrac{2^{n+2}}{3}\pi$. Similarly, if $k\leq -1$,
then $2^j(\xi + 2k\pi)\leq -\tfrac{2^{n+2}}{3}\pi$. Hence, for
$j\geq n+1$, only contributing $k$ to $D_{\psi}$ is $k=0$.

Thus, we have
 \begin{eqnarray}
 D_{\psi}(\xi ) & = &
 \sum_{j=1}^{n-1}\sum_{k\in\Z}|\hat\psi
 \left(2^j(\xi+2k\pi)\right)|^2 \nonumber \\
 & & +\sum_{k=-1}^1|\hat\psi\bigl(2^n(\xi + 2k\pi)\bigr)|^2
 +\sum_{j\geq n+1}|\hat\psi(2^j\xi)|^2 \label{eqn:dpsi2} \\
 & = &
 \sum_{j=1}^{n-1}\sum_{k\in\Z}|\hat\psi
 \left(2^j(\xi+2k\pi)\right)|^2 \nonumber \\
 & & + \left\{|\hat\psi \left(2^n(\xi-2\pi)\right)|^2
 +|\hat\psi \left(2^n(\xi + 2\pi)\right)|^2\right\} 
 +\sum_{j\geq n}|\hat\psi(2^j\xi)|^2. \label{eqn:dpsi1}
\end{eqnarray}

If $ \xi\in [-\tfrac{2}{3}\pi,\tfrac{2}{3}\pi]$, then
 $2^n(\xi+2\pi)\geq 2^n\cdot \tfrac{4}{3}\pi
 =\tfrac{2^{n+2}}{3}\pi$.
Similarly, $2^n(\xi - 2\pi)\leq -\tfrac{2^{n+2}}{3}\pi$. So both
the terms inside the curly bracket in (\ref{eqn:dpsi1}) are zero, and we get $D_{\psi}\leq (n-1)+1=n$. Now, if $\xi\in [\tfrac{2}{3}\pi,\pi]$, then for all $j\geq n+1$, we have 
$2^j\xi\geq 2^{n+1}\cdot \tfrac{2}{3}\pi = \tfrac{2^{n+2}}{3}\pi$.
So the last sum in (\ref{eqn:dpsi2}) is zero and again
$D_{\psi}\leq n$. In a similar manner, it can be shown that
$D_{\psi}\leq n$ if $\xi\in [-\pi,-\tfrac{2}{3}\pi]$. This finishes the proof.
\qed

\proof[Proof of Theorem {\rm\ref{thm:main}}]
The wavelets we construct to prove Theorem \ref{thm:main} are MSF wavelets so that it suffices to construct the associated wavelet sets. In addition to proving the theorem, this construction also provides an example of a family of wavelet sets which are union of infinite number of intervals. We will treat the even and odd cases separately.

\vskip 2mm
\noindent{\bf Case I}. $n$ is even

\medskip

 For a real number $\ep$ such that
 $0< \ep < \delta=\tfrac{2^{n+2}}{3(2^{n+2}-1)}\pi$, let
 $S_i,\ 1\leq i\leq 6$ be the following sets.
\begin{eqnarray*}
    S_1 &=& \bigl[-\tn, -\tn+\ep\bigr], \\
    S_2 &=& \bigl[-\tfrac{1}{3}\pi+\two, -\tfrac{1}{6}\pi\bigr], \\
    S_3 &=& \bigl[\tn+\ep-2\pi, \tn-\tfrac{5}{3}\pi+\two\bigr], \\
    S_4 &=& \bigl[\tn-\tfrac{3}{2}\pi, \tn-\tfrac{7}{6}\pi+\thr\bigr], \\
    S_5 &=& \bigl[\tn-\pi, \tn-\tfrac{2}{3}\pi\bigr], \\
    S_6 &=& \bigl[\tn-\tfrac{2}{3}\pi+\ep, \tn+\ep\bigr].
\end{eqnarray*}
Define the sets $X_0$, $Y_0$ and $Z_0$ as
\begin{eqnarray*}
X_0 &=& \bigl[\tfrac{1}{6}\pi+\thr,            
         \tfrac{1}{3}\pi-\tfrac{1}{2^{n+1}}\pi+\two\bigr], \\
Y_0 &=& \tfrac{1}{2^{n+2}}\bigl(S_2+2\cdot\tfrac{2^{n+1}-2}{3}\pi\bigr), \\
Z_0 &=& \bigl[\tfrac{1}{3}\pi-\tfrac{1}{3\cdot 2^{n+1}}\pi,
         \tfrac{1}{3}\pi-\tfrac{1}{3\cdot 2^{n+1}}\pi+\two\bigr].
\end{eqnarray*}
The parameter $\ep$ is chosen in a suitable manner to make the above sets non-empty. For $j\geq 1$, let the sets $X_j,Y_j,Z_j$ be defined recursively
as follows:
\[
 P_j = \tfrac{1}{2^{n+2}}\bigl(P_{j-1}+2\cdot\tfrac{2^{n+1}-2}{3}\pi\bigr),\  j\geq 1, \quad P\in\{X,Y,Z\}.
\]
By a routine calculation, we can easily verify the following facts:
\begin{enumerate}
\item[(i)] 
  $P_j\subset\bigl[\tfrac{1}{6}\pi+\thr, \tfrac{1}{3} \pi\bigr], j\geq 0,
  \quad P\in\{X,Y,Z\}$.
\item[(ii)] $2^{n+2}P_j\subset\bigl[\tn-\tfrac{7}{6}\pi+\thr,   
  \tn-\pi\bigr]$ for $j\geq 1, \quad P\in\{X,Y,Z\}$.
\item[(iii)] $\{X_j,Y_j,Z_j:j\geq 0\}$ is a disjoint collection.
\item[(iv)] $X_j$ lies to the left of $Y_j$, and $Y_j$ lies to the left
  of $Z_j$ for $j\geq 0$.
\item[(v)] $X_{j+1},Y_{j+1},Z_{j+1}$ lie between $Y_j$ and $Z_j,
  j\geq 0$.
\end{enumerate}

 Let
\[
 X=\bigcup_{j\geq 0}X_j,\ Y=\bigcup_{j\geq 0}Y_j,\
 Z=\bigcup_{j\geq 0}Z_j,
\]
and
\[
 V=[\tn-\tfrac{7}{6}\pi+\thr, \tn-\pi]\setminus
 \Bigl\{\bigcup_{j\geq 1}2^{n+2}(X_j\cup Y_j\cup Z_j)\Bigr\}.
\]
Now define
\begin{equation}\label{eq:wavset}
 W=\Bigl(\bigcup_{i=1}^{6}S_i\Bigr)\cup(X\cup Y\cup Z)\cup V.
\end{equation}

 \noindent{\bf Claim.} $W$ is a wavelet set.

{\it Translation equivalence}: 
We will show that $W$ is translation equivalent to the interval $[\tn+\ep-2\pi,\tn+\ep]$ of length $2\pi$. Note that
\[
 S_3\cup \bigl(S_2+2\cdot\tfrac{2^{n+1}-2}{3}\pi\bigr)\cup S_4=
 \bigl[\tn+\ep-2\pi,\tn-\tfrac{7}{6}\pi+\thr\bigr],
\] and
\[
 S_5\cup \bigl(S_1+2\cdot\tfrac{2^{n+2}-1}{3}\pi\bigr)\cup S_6=
 \bigl[\tn-\pi,\tn+\ep\bigr].
\]
Now,
\begin{eqnarray*}
 \lefteqn{(X\cup Y\cup Z)+2\cdot\tfrac{2^{n+2}-1}{3}\pi} \\
&=& \bigcup_{j\geq 0} 
    \Bigl\{(X_j+2\cdot\tfrac{2^{n+1}-2}{3}\pi)\cup
  (Y_j+2\cdot\tfrac{2^{n+1}-2}{3}\pi)\cup
  (Z_j+2\cdot\tfrac{2^{n+1}-2}{3}\pi)\Bigr\} \\
& = & \bigcup_{j\geq 0}2^{n+2}( X_{j+1}\cup Y_{j+1}\cup Z_{j+1}) 
  = \bigcup_{j\geq 1}2^{n+2}( X_j\cup Y_j\cup Z_j).
\end{eqnarray*}
 Therefore, by the definition of $V$
 \[ V\cup\bigl\{(X\cup Y\cup Z)+2\cdot\tfrac{2^{n+1}-2}{3}\pi\bigr\}=
 \bigl[\tn-\tfrac{7}{6}\pi+\thr,\tn-\pi\bigr]. \]
 
Observe that $\tfrac{2^{n+1}-2}{3}\pi$ and $\tfrac{2^{n+2}-1}{3}\pi$ are
integers, since $n$ is even. We have proved that appropriate translations of the partition of $W$ in (\ref{eq:wavset}) form a partition of the interval $[\tn+\ep-2\pi,\tn+\ep]$. So $W$ is translation equivalent to this interval. 

\vskip 2mm

{\it Dilation equivalence}:
It is enough to show that $W\cap(-\infty,0)$ and
 $W\cap(0,\infty)$ are respectively dilation equivalent to the
 intervals $[-\tn,-\tfrac{2^{n+1}}{3}\pi]$ and
$[\tfrac{2^{n+1}}{3}\pi+\tfrac{\ep}{2},\tn+\ep]$.
\begin{eqnarray*}
 S_1\cup (2^{n+2}S_2) &=& \bigl[-\tn,-\tfrac{2^{n+1}}{3}\pi\bigr], \\
(2^{n+2}X_0)\cup S_3 \cup(2^{n+2}Y_0)\cup S_4
 & = & \bigl[\tfrac{2^{n+1}}{3}\pi+\tfrac{\ep}{2},
 \tn-\tfrac{7}{6}\pi+\thr\bigr], \\
S_5 \cup (2^{n+2}Z_0)\cup S_6 & = & \bigl[\tn-\pi, \tn+\ep\bigr], \\
 \Bigl\{2^{n+2}\Bigl(\bigcup_{j\geq 1}(X_j\cup Y_j\cup   
  Z_j)\Bigr)\Bigr\}\cup V
& = & \bigl[\tn-\tfrac{7}{6}\pi+\thr,\tn-\pi\bigr].
\end{eqnarray*}
Hence, $W$ is a wavelet set.

\medskip
 \noindent{\bf Case II}. $n$ is odd
\medskip

This case is dealt in a similar manner, but we have to start with
different sets. For $0< \ep < \tfrac{2^{n+2}}{3(2^{n+2}-1)}\pi$,
let the sets $S_1$, $S_2$ be as above. Define
\begin{eqnarray*}
S_3 &=& \bigl[\tn+\ep-2\pi, \tn-\tfrac{4}{3}\pi\bigr], \\
S_4 &=& \bigl[\tn-\tfrac{4}{3}\pi+\ep, \tn-\pi+\two\bigr], \\
S_5 &=& \bigl[\tn-\tfrac{5}{6}\pi, \tn-\tfrac{1}{2}\pi+\thr\bigr], \\
S_6 &=& \bigl[\tn-\tfrac{1}{3}\pi, \tn+\ep\bigr].
\end{eqnarray*}
Let
\begin{eqnarray*}
X_0 &=& \bigl[\tfrac{1}{6}\pi+\thr,   
  \tfrac{1}{3}\pi-\tfrac{1}{2^{n+1}}\pi+\two\bigr], \\
Y_0 &=& \bigl[\tfrac{1}{3}\pi-\tfrac{1}{3\cdot 2^n}\pi,
  \tfrac{1}{3}\pi-\tfrac{1}{3\cdot 2^n}\pi+\two\bigr], \\
Z_0 &=& \tfrac{1}{2^{n+2}}\bigl(S_2+2\cdot\tfrac{2^{n+1}-1}{3}\pi\bigr).
\end{eqnarray*}
In this case also, the choice of $\ep$ ensures that the above sets are non-empty. Define the sets $X,\ Y$ and $Z$ as in Case I. Let
\[
 V=\bigl[\tn+\thr-\tfrac{1}{2}\pi,\tn-\tfrac{1}{3}\pi\bigr]\setminus
 \Bigl\{\bigcup_{j\geq 1}2^{n+2}\bigl(X_j\cup Y_j\cup Z_j\bigr)\Bigr\},
\]
and let $W$ be defined by (\ref{eq:wavset}).

As in the case when $n$ is even, to show that $W$ is a wavelet set, we show the translation equivalence of $W$ with the interval
$[\tn+\ep-2\pi,\tn+\ep]$; and the dilation  equivalence of
$W\cap(-\infty,0)$ and $W\cap(0,\infty)$ with the intervals
$[-\tn,-\tfrac{2^{n+1}}{3}\pi]$ and
$[\tfrac{2^{n+1}}{3}\pi+\tfrac{\ep}{2},\tn+\ep]$ respectively.

To see the translation equivalence, observe that
\[
 S_3\cup \bigl(S_1+2\cdot\tfrac{2^{n+2}-2}{3}\pi\bigr)\cup S_4
 \cup \bigl(S_2+2\cdot\tfrac{2^{n+1}-1}{3}\pi\bigr)\cup S_5=
\]
\[
 \bigl[\tn+\ep-2\pi,\tn+\ep\bigr],
\]
\[
S_6= [\tn-\tfrac{1}{3}\pi, \tn+\ep].
\]
It can be shown, in a manner similar to Case I, that
\[
 V\cup\bigl\{(X\cup Y\cup Z)+2\cdot\tfrac{2^{n+1}-1}{3}\pi\bigr\}=
 \bigl[\tn-\tfrac{1}{2}\pi+\thr,\tn-\tfrac{1}{3}\pi\bigr].
\]

For dilation equivalence, we observe
\[
 S_1\cup (2^{n+2}S_2) = \bigl[-\tn,-\tfrac{2^{n+1}}{3}\pi\bigr], \]
\[ (2^{n+2}X_0)\cup S_3 \cup(2^{n+2}Y_0)\cup S_4
  \cup (2^{n+2}Z_0)\cup S_5 = \]
\[ \bigl[\tfrac{2^{n+1}}{3}\pi+\tfrac{\ep}{2},
  \tn-\tfrac{1}{2}\pi+\thr\bigr], \]
\[  S_6 = \bigl[\tn-\tfrac{1}{3}\pi, \tn+\ep\bigr], \]
\[ \Bigl\{2^{n+2}\Bigl(\bigcup_{j\geq 1}(X_j\cup Y_j\cup Z_j)\Bigr)\Bigr\}\cup V
 = \bigl[\tn-\tfrac{1}{2}\pi+\thr,\tn-\tfrac{1}{3}\pi\bigr].\]
Hence, in this case also we have proved that $W$ is a wavelet set.

By definining $\hat\psi=\chi_W$, we get a wavelet $\psi$ such that $\hat\psi$ is supported in $[-\tn,\tn+\ep]$, since $W$ is a subset of this interval.

Finally, to complete the proof of Theorem \ref{thm:main}, we have to show that $\|D_\psi\|_\infty>n$, where
$\hat\psi=\chi_W$. We prove $D_\psi(\xi)\geq n+1$ for  a.e.\
$\xi\in[\tfrac{2}{3}\pi,\tfrac{2}{3}\pi+\tfrac{\ep}{2^{n+1}}]$.

For $1\leq j\leq n+1$, let $k_j=\tfrac{2^{(n+1-j)}-1}{3}$ and
$l_j=-\tfrac{2^{(n+1-j)}+1}{3}$. Observe that $k_j$ is an integer
if $n-j$ is odd, and $l_j$ is an integer when $n-j$ is even.

Let
$\xi\in[\tfrac{2}{3}\pi,\tfrac{2}{3}\pi+\tfrac{\ep}{2^{n+1}}]$. 
If $n$ is even, then for $j=1,3,\dots,n+1$, we have
$2^j(\xi+2k_j\pi)=2^j(\xi-\frac{2}{3}\pi)+
\tfrac{2^{n+2}}{3}\pi\in[\tn, \tn+\ep]\subset S_6$.
Also, for $j=2,4,6,\dots n$, $2^j(\xi+2l_j\pi)\in [-\tn,-\tn+\ep]=S_1$. 
Similarly, if $n$ is odd, then
$2^j(\xi+2k_j\pi)\in[\tn, \tn+\ep]$ if $j=2,4,6,\dots n+1$, and
$2^j(\xi+2l_j\pi)\in[-\tn, -\tn+\ep]$ if $j=1,3,5,\dots,n$.

In each case, there are $n+1$ different pairs of $(j,k)$, with
$j\geq 1$ and $k\in\Z$, such that $2^j(\xi+2k\pi)$ is contained in $W$ which is
the support of $\hat\psi$. Each such pair will contribute $1$ to
the sum $D_\psi(\xi)$ defined in (\ref{dpsi}). Therefore, $\|D_\psi\|_\infty\geq n+1$.
\qed
%%===================================

\end{document}